\documentclass[12pt]{article}
\usepackage{amssymb,amsmath,amsfonts,amsthm}

 \newtheorem{theorem}{Theorem}[section]

\begin{document}

\title{Nonexistence results for a class of nonlinear elliptic equations involving critical Sobolev exponents  }

\author{Cristina Tarsi \thanks{ e-mail: tarsi@mat.unimi.it. \ The
author is member of the research group G.N.A.M.P.A of the Italian
Istituto Nazionale di Alta Matematica (INdAM).} ,\\
Dipartimento di Matematica, Universit\`{a} degli Studi, \\
I-20133 Milano, Italy}
\date{}

\maketitle


% main text
\section{Introduction}

Let $\Omega $ be a bounded domain in $\mathbb{R}^N$, with $N\geq
3$;
 consider the following semilinear elliptic problem
\begin{equation}
\label{1}
\left\{
\begin{array}{ll}
 -\Delta u=\lambda g\left( x,u\right) +\left| u\right| ^{2^{*}-2}u,%
& \;\;\;\;\;x\in \Omega \\
 u>0 & \;\;\;\;\;x\in\Omega\\
 u=0 & \;\;\;\;\;x\in \partial \Omega
\end{array}
\right.
\end{equation}
where $2^{*}=2N/\left( N-2\right) $ is critical from the viewpoint of the
Sobolev embedding $H_0^1\left( \Omega \right) \subset L^{2^{*}}\left( \Omega
\right) $, and $g\left( x,u\right) $ is a lower-order perturbation of $%
u^{2^{*}-1}$, in the sense that $\lim_{u\rightarrow +\infty
}g\left( x,u\right) /u^{2^{*}-1}=0$. As well known, if $g$
satisfies suitable assumptions, solutions of (\ref{1})
 correspond
to critical points of the functional
\[
\Psi \left( u\right) =\frac 12\int_\Omega \left| \nabla u\right|
^2dx-\lambda \int_\Omega G\left( x,u\right) dx-\frac 1{2^{*}}\int_\Omega
\left| u\right| ^{2^{*}}dx,
\]
where $G\left( x,u\right) =\int_0^ug\left( x,t\right) dt$. Since
the embedding $H_0^1\left( \Omega \right) \subset L^{2^{*}}\left(
\Omega \right) $ is not compact, the functional $\Psi $ does not
satisfy the Palais-Smale condition: hence the standard variational
arguments do not apply. For equations with critical growth,
nontrivial solution may non exists: a
well-known nonexistence result due to Pohozaev \cite{P} asserts that if $%
\Omega $ is starshaped and $\lambda \leq 0$ there is no solution (different
from the trivial one) of the problem
\begin{equation}
\left\{
\begin{array}{ll}
-\Delta u=\lambda u+\left| u\right| ^{2^{*}-2}u & \;\;\;\;\;x\in%
\Omega \\
u>0 & \;\;\;\;\;x\in\Omega \\
u=0 & \;\;\;\;\;x\in\partial \Omega .
\end{array}
\right.
\end{equation}
In recent years this situation of lack of compactness has been
extensively investigated (see for example \cite{N}); according to
the behaviour of $g$ and the kind of results one wants to prove,
topological or variational methods turn out to be more
appropriate. When $g$ is superlinear, for example $g=\left|
u\right| ^{p-1}u$, $1<p<2^{*}-1$, variational tools, such as
minimax arguments, provide the existence of a nontrivial positive
solution; on the contrary when $g$ is sublinear, for example
$g=\left| u\right| ^{p-1}u$, $0<p<1$, sub- and super-solutions are
quite convenient.
In particular we recall the following known existence results for problem $%
\left( \ref{1}\right) $:

\begin{itemize}
\item  The first existence result is due to Brezis-Nirenberg \cite{BN}; in
a pioneering result, they showed that, when $g\left( x,u\right)
=u$, there exists a nontrivial positive solution if $\lambda \in
\left( \lambda ^{*},\lambda _1\right) $, with $\lambda ^{*}=0$ for
$N\geq 4$ and $0<\lambda ^{*}\left( \Omega \right) <\lambda _1$
for $N=3$ ($\lambda _1$ denoting the first eigenvalue of $-\Delta
$ relative to the homogeneous Dirichlet problem in $\Omega $). In
the same work they also proved an existence result for equation
(\ref{1}) when $g$, roughly speaking, has a linear or superlinear
growth near zero and near infinity: in this case there is again
bifurcation from infinity in $\lambda =0$ for $N\geq 4$ , whereas
for $N=3$ it can not be guaranteed in the entire subcritical
growth range of the term $g$.

\item  Later, Ambrosetti-Brezis-Cerami \cite{ABC} established the
existence of two positive solution for $0<\lambda <\Lambda $ when $g=u^q$
with $0<q<1$ and $N\geq 3$, thanks to the combined effects of the sublinear
and superlinear terms. The first solution is found using sub- and
super-solutions; in contrast with the pure concave case, a second positive
solution is found by variational arguments. Moreover, they proved that the
first solution, $u_\lambda $, is such that $\left\| u_\lambda \right\|
_\infty \rightarrow 0$ as $\lambda \downarrow 0$, while the second solution,
$w_\lambda $, (if $\Omega $ is strictly starshaped) has a nonlimited norm,
that is, $\left\| w_\lambda \right\| _\infty \rightarrow 0$ as $\lambda
\downarrow 0.$
\end{itemize}

One may ask if the superlinear/sublinear growth of the subcritical
term can be weakened in these existence results, e. g.,
considering subcritical terms presenting superlinear or subliner
asymptotic behaviour near zero or near infinity. In fact, the
proofs presented by Brezis-Nirenberg in \cite{BN} and by
Ambrosetti-Brezis-Cerami \cite{ABC} can be generalized with some
technicalities to subcritical terms presenting, respectively, a
superlinear or sublinear asymptotic behaviour near the origin: for
example, it is not hard to prove the existence of a positive
solution for problem (\ref{1}) if $N\geq 5$ and $g(x,u)$ satisfies
the following assumptions

\begin{equation}\label{c1ip3}
\left\{
\begin{array}{ll}
 (i)\;\;\; g\left( x,u\right) =\left| u\right| ^{p-1}u & \;\;\;\textrm{for}\;\;\left|
u\right| <1,\;\;%
x\in \Omega ,\;\;\; p>1 \\
 & \\
(ii)\;\;\; \exists \delta >0:g\left( x,s\right) \geq 0 &
\;\;\;\;\; \forall
\left| x\right| <\delta ,\;\;\;\forall s>0,\\
\end{array}
\right.
\end{equation}

and the existence of two positive solutions for
$\lambda\in(0,\Lambda) $ if $N\geq 3$ and $g(x,u)$ satisfies

\begin{equation}
\left\{
\begin{array}{ll}
(i)\;\;\; g\left( x,u\right) =\left| u\right| ^{p-1}u &
\;\;\;\textrm{for}\;\;\left| u\right| <1,\;\;%
x\in \Omega, \\
&  0<p<1,\;\;\;p<p_1<2^{*}-1
\\
(ii)\left| u\right| ^{p-1}u\leq g\left( x,u\right) \leq \left|
u\right| ^{p_1-1}u & \;\;\;\textrm{for}\;\; \left| u\right| \geq
1,\;\;\;x\in \Omega .
\end{array}
\right. \label{c1ip4}
\end{equation}

We note that the behaviour of the subcritical term near the origin
seems to determine the structure of bifurcation from infinity for
problem (\ref{1}): that is, it seems not possible to obtain
similar existence results assuming only superlinear or sublinear
growth near infinity. In lower dimensions, however, the effect of
the pure convex/concave behaviour of the subcritical term assumes
an increasing role, that can not be replaced by the analogous
asymptotic behaviour of $g$ near zero: for example, in the pure
convex case considered in \cite{BN} the existence results are
valid for $N\geq4$, whereas for $N=3$ a nonexistence result is
given; assuming convexity near the origin, instead, the existence
results can be extended, in general, only for $N\geq5$, as we will
prove exhibiting a counterexample. The aim of this paper is to
point out the difference between the pure convex/concave case and
the case of convex/concave growth of the subcritical term near the
origin in lower dimensions: based on a celebrated Identity due to
Pohozaev \cite{P}, we construct special classes of nonlinear
problems which do not have nontrivial solutions bifurcating from
infinity in $\lambda =0$ (if the domain $\Omega $ is strictly
starshaped), according to the behaviour of the subcritical term
$g\left( x,u\right) $ and to the dimension $N$. In particular, we
prove that the first critical dimension is $N=4$, which is somehow
in contrast with the pure convex case considered in \cite{BN}. We
remark that the class of subcritical terms presented here has
superlinear growth near the origin, whereas the growth near
infinity can be sublinear or superlinear; sublinear growth near
the origin, instead, determines bifurcation from infinity for all
$N\geq3$, either for convex or for concave behaviour of $g$ near
infinity, as one can prove following \cite{ABC} with slight
modifications: that is, the role of the asymptotic behaviour of
the subcritical term $g(x,u)$ near infinity does not determine the
structure of bifurcation of problem (\ref{1}).

\section{Recalls from potential theory and elliptic estimates}

Let $\Omega $ be a bounded (smooth) domain in $\mathbb{R} ^N$,
with $N\geq 3$. We will exhibit two classes of subcritical terms
$g(x,u)$ such that problem (\ref{1}) does not admit any positive
solution when $\lambda$ is close to zero and $N=3,4$. The proofs
rely on the so-called \emph{Pohozaev's identity }\cite{P}: suppose
$u$ is a smooth function satisfying
\begin{equation}
\left\{
\begin{array}{ll}
 -\Delta u=f(u)%
& \;\;\;\;\;x\in \Omega \\
  u=0 & \;\;\;\;\;x\in \partial \Omega
\end{array}
\right.
\end{equation}
where $g$ is a continuous function on $\mathbb{R}$ and $\Omega $
is a (smooth) starshaped domain. Then we have
\begin{equation}
\label{pohoz}
\left( 1-\frac 12n\right) \int_\Omega f\left(
u\right) \cdot udx+n\int_\Omega F\left( u\right) dx=\frac
12\int_\Omega \left( x\cdot \nu \right) \left( \frac{\partial
u}{\partial \nu }\right) ^2ds
\end{equation}
where
\[
F\left( u\right) =\int_0^uf\left( t\right) dt
\]
and $\nu $ denotes the outward normal to $\partial \Omega .$ We
will combine the Pohozaev's identity (\ref{pohoz}) together with
some standard elliptic inequalities and the weak interpolation
inequality, which we briefly recall in the following (see
\cite{LL}).
 \\
Let us consider a domain $\Omega\subseteq\ \mathbb{R} ^N$; denote
with $\emph{D}(\Omega)$ the space of the test functions and with
$\emph{D'}(\Omega)$ the space of distributions, that is, the dual
space of $\emph{D}(\Omega)$. Let us recall the definition of the
Green's functions for the Poisson's equation in $\mathbb{R} ^N$,
\begin{eqnarray}
\label{poisson}
G(x)&=&-\frac{1}{2\pi}\ln|x|
\;\;\;\;\;\;\;\;\;\;\;\;\;\;\;\;\;\;\;\;\;N=2\nonumber
\\
g(x)&=&\frac{1}{(N-2)\mu(\mathbb{S}^{N-1})}|x|^{2-N}\;\;N\neq2
\end{eqnarray}
where $\mu(\mathbb{S}^{N-1})$ is the area of the unit sphere
$\mathbb{S}^{N-1}\subseteq \mathbb{R} ^N$. It is well known that
for every $u\in L^{1}_{\mathrm{loc}}$, the function

\[
k_{u}(x)=(G\ast u)(x)=\int_{\Omega}G(x-y)u(y)dy
\]
satisfies

\begin{eqnarray}
\label{2}
 k_{u}&\in & L^{1}_{\mathrm{loc}}(\Omega) \nonumber
\\
\\ \nonumber
-\triangle k_{u}&=& u \;\;\; \in \emph{D'}(\Omega)
\end{eqnarray}

if the function $y\mapsto G(x-y)u(y)$ is summable over $\Omega$
for almost every $x$. On the other hand, applying the Young's
inequality

\[
\| g\ast h\| _{p}\leq C_{q,r,p,N}\| g\| _{q}\| h\| _{r}
\;\;\;\mathrm{if}\;\;
\frac{1}{q}+\frac{1}{r}=1+\frac{1}{p},\;\;\;p,q,r\geq1
\]

with $g=G,h=u$ and $r=1$, we have that

\begin{equation}
\label{3}
 \| k_{u}\| _{p} \leq C_{q,r,p,N} \| G\| _{p} \| u\|
_{1}.
\end{equation}

Therefore, combining (\ref{poisson}),(\ref{2}) and (\ref{3}) we
can conclude that the operator $\triangle ^{-1}$ is bounded from
$L^{1}$ to $L^{p}$ with $p\in [1,3)$ if $N=3$, and from $L^{1}$ to
$L^{p}$ with $p\in [1,2)$ if $N=4$; that is, for every $v \in
L^{1}$ there is $u \in L^{p}$ (with $p$ satisfying the previous
conditions) such that

\begin{eqnarray}
\label{4}
 \triangle u &=& v \nonumber
\\
\\ \nonumber
\| u\| _{p} &\leq & C_{p,N} \| v\| _{1} =C_{p,N}\| \triangle
u\|_{1}.
\end{eqnarray}

If $p=3$ and $N=3$, or, respectively, if $p=2$ and $N=4$,
(\ref{4}) are not verified; in this case, however, we can apply
the notion of weak $L^{p}$ spaces (see \cite{LL}). Consider the
space of all measurable functions $u$ such that

\begin{equation}\label{wLp}
[u]_{q,w} =\sup_{\alpha >0} \;\;\alpha \cdot \mu \{x:|u(x)|>\alpha
\}^{1/q}<\infty;
\end{equation}

this space is called weak $L^{q}$-space $L^{q}_{w}(\mathbb{R}
^{N})$. Any function in $L^{q}(\mathbb{R} ^{N})$ is in
$L^{q}_{w}(\mathbb{R} ^{N})$: simply note that

\[
\| u\|^{q}_{q}\geq \int_{|u|>\alpha }|u(x)|^{q}dx\geq \alpha
^{q}\cdot \mu \{x:|u(x)|>\alpha \}=[u]_{q,w}^q.
\]

The expression (\ref{wLp}) does not define a norm; nevertheless,
there is an alternative expression, equivalent to (\ref{wLp}),
that is indeed a norm: it is given by

\begin{equation}\label{wnorm}
\| u\| _{q,w}=\sup_{A}\frac{1}{\mu (A) ^{1/r}}\int_A| u(x)|dx.
\end{equation}

where $1/q+1/r=1$ and $A$ denotes an arbitrary measurable set of
measure $\mu (A)<\infty$. In particular, $u(x)=|x|^{-\lambda}$ is
in $L^{q}_{w} (\mathbb{R} ^N)$ with $q=N/\lambda$, $N>\lambda >0$
and
\[
\| u\| _{N/\lambda ,w}=\frac{N}{N-\lambda}\left[\frac{\mu (\mathbb
S^{N-1})}{N}\right]^{\lambda /N}.
\]
The weak Young inequality states that for $g \in
L^{q}_{w}(\mathbb{R} ^N)$ and $\infty>p,q,r>1$ with
$1/p+1/q+1/r=2$, the following inequality holds:

\begin{equation}\label{young}
\int _{\mathbb{R} ^N} \int _{\mathbb{R} ^N} f(x)g(x-y)h(y)dxdy\leq
C_{p,q,r}\| f\| _{p}\| g\| _{q,w}\| h\| _{r};
\end{equation}

taking $\lambda =N/q$ and $g(x)=|x|^{-\lambda}$ the weak Young
inequality (\ref{young}) is equivalent to the
Hardy-Littlewod-Sobolev inequality,

\[
\int _{\mathbb{R} ^N} \int _{\mathbb{R} ^N}
f(x)|x-y|^{-\lambda}h(y)dxdy\leq C_{N,\lambda ,p}\| f\| _{p}\| h\|
_{r}
\]

with $p,r>1,0<\lambda <N$ and $1/p+\lambda /N+1/r=2$. In
particular, the sharp constant in the weak Young inequality is the
same as for the Hardy-Littlewood-Sobolev inequality. Observe that
we can also view the Young inequality as the statement that the
convolution is a bounded map from $L^{p}(\mathbb{R} ^N)\times
L^{q}_{w} (\mathbb{R} ^N)$ to $L^{s}(\mathbb{R} ^N)$, that is

\begin{equation}\label{conv}
\| f\ast g\|_{s}\leq C_{q,s,p,N}\| g\| _{q,w} \| f\| _{p}\;\;\;
\mathrm{if} \;\;\; \frac{1}{p}+\frac{1}{q}=1+\frac{1}{s},\;\;\;
p,q,r>1.
\end{equation}

A final inequality involving the $L^{q}_{w}$ spaces is the weak
interpolation inequality: if $u\in L^{p}_{w}\cap L^{r}$, with
$r<p$, then

\begin{equation}\label{winterp}
\| u\| _q\leq K_{q,r,p,N}\| u\| _{p,w}^a\| u\|
_r^{1-a}\hspace{0.2in}\textrm{with } \frac 1q=\frac
ap+\frac{1-a}r.
\end{equation}

This inequality will allow us to combine the estimates obtained
from the Pohozaev's identity with the elliptic estimates
(\ref{4}).

\section{Nonexistence results}

In this section we construct two classes of nonlinear elliptic
problems with critical growth which don't admit any positive
solution near $\lambda =0$. The proof of nonexistence is based on
Pohozaev's identity and on the elliptic estimates presented in the
previous section. From now on suppose $\Omega$ is strictly
starshaped about the origin, so that $\left( x\cdot \nu \right)
>c>0$ a.e. on $\partial\Omega .$ We discuss separately the two
cases, $N=3$ and $N=4.$

\subsection{The case $N=3$.}

We assume here that $N=3$ and
\begin{equation}\label{N=3}
g( u) =\left\{
\begin{array}{l}
\| u\| ^{p-1}\cdot u\hspace{0.5in}| u| <1, \;\;1<p \\
\| u\| ^{q-1}\cdot u\hspace{0.5in}| u| \geq 1,\;\; 0<q\leq 3
\end{array}
\right.
\end{equation}
Then we have the following result.

\begin{theorem}\label{t1}
Let $\Omega $ be strictly starshaped about the origin; suppose
that $u$ is a solution of problem (\ref{1}), with $g$ given by
(\ref{N=3}). Then
\[
\lambda \geq \lambda _0\left( q,p,\Omega \right) >0
\]
if $1<p,$ $0<q\leq 3$.
\end{theorem}

\emph{Proof of Theorem \ref{t1}. }By Pohozaev's identity
 (\ref{pohoz}), since $\Omega $ is strictly starshaped, we have
\begin{equation}\label{5}
\begin{array}{l}
\lambda \displaystyle \frac{5-q}{2(q+1)}\int_{|u| \geq 1}|u|
^{q+1}dx+\lambda
\frac{5-p}{2(p+1)}\int_{|u| <1}|u| ^{p+1}dx \\
\\
+3\lambda \cdot \mu \left\{x\in \Omega \: |u| \geq
1\right\} \displaystyle \frac{q-p}{(q+1)(p+1)}=\frac 12\int_{\partial \Omega}%
 (x,\nu ) \left| \frac{\partial u}{\partial \nu }\right| ^2dx\geq \\
 \\
c\left( \int_\Omega \left| \Delta u\right| dx\right) ^2.
\end{array}
\end{equation}
We discuss separately the different cases.

\begin{enumerate}
\item[(i)]
If $1<q=p\leq 3$, the subcritical term $g$ defined by (\ref{N=3})
reduces to the case considered in Theroem $2.4$ in \cite{BN}, so
we will be brief. Indeed, combining (\ref{5}) with equation
(\ref{1}) implies

\begin{equation}\label{6}
\lambda c(q) \| u\|^{q+1}_{q+1}\geq c\| \triangle u\|^{2}_{1}\geq
\| u\|^{10}_{5};
\end{equation}

on the other hand, by the elliptic estimates (\ref{4}) we obtain

\[
\| \triangle u\|^{2}_{1}\geq c\| u\|^{2}_{r}
\]

for all $r\in [1,3)$, so that

\[
\lambda  \| u\|^{q+1}_{q+1}\geq c\| u\|^{2}_{r}.
\]

Using the interpolation inequality

\[
\| u\|_{q+1}\leq \| u\|^{a}_{r}\cdot \| u\|^{1-a}_{5}
\]

with $\displaystyle \frac 1{q+1} = \frac ar + \frac{1-a}5$,
$r<q+1<5$, and combining with the previous inequality, one finds
\[
\| u\|_{q+1}\leq c\lambda ^{\frac{4a+1}{10}}\|
u\|^{(q+1)(4a+1)/10}_{q+1}.
\]
Choosing
\[
r=\frac{9-q}{5-q},\;\; a=\frac{9-q}{4(q+1)}
\]
which satisfy $a<1,r<q+1,r\in [1,3)$ if $1<q<3$ we obtain
\[
\lambda ^{1/(q+1)}\geq c.
\]

\item[(ii)]
 If $q=p=3$, we obtain $r=3$, so that inequality $\|
\triangle u\|^{2}_{1}\geq c\| u\|^{2}_{r}$ is not verified. In
this case, we can observe that, by the potential theory,
\[
u\leq v=\frac{c}{|x|}\ast |\triangle u|,
\]

and $|x|^{-1}\in L^3_{w}$; combining these two relations with
(\ref{wnorm}) yields the following inequality

\begin{equation}\label{7}
\| \triangle u\|^{2}_{1}\geq c\| u\|^{2}_{3,w}.
\end{equation}

 Then, the proof can be completed as for $q=p<3$ using the
 following weak interpolation inequality
\[
\| u\|_{4}\leq \| u\|^{3/8}_{3,w}\cdot c\| u\|^{5/8}_{5}.
\]

\item[(iii)]
If $3\geq q>p>1$, then $q-p>0$ and

\begin{equation}\label{7a}
\| u\| _{q+1}^{q+1}\geq \int_{|u| \geq 1}|u| ^{q+1}dx>\mu \left\{
 x\in \Omega :|u| \geq 1\right\} .
\end{equation}

On the other hand,

\begin{equation}\label{8}
\int_{|u|<1}|u| ^{p+1}dx=\int_{|u|<1}|u| ^{q+1}dx+\int_{|u|
<1}\left( |u| ^{p+1}-|u|^{q+1}\right) dx
\end{equation}

Let us estimate the second integral in the right hand side of
(\ref{8}) as follows:

\begin{eqnarray*}
0 &<&\int_{|u| <1}\left( |u|^{p+1}-|u|^{q+1}\right) dx=\int_\Omega
|u| \cdot \chi _{\left\{ |u| <1\right\} }\cdot \left( |u|^p-|u|^q\right) dx \\
&\leq &\left\| u\right\| _{q+1}\cdot \left\{ \int_{|u|<1}\left[
|u| ^q\cdot \left( |u| ^{p-q}+1\right) \right] ^{%
\frac{q+1}q}dx\right\} ^{\frac q{q+1}} \\
&\leq& 2\left\| u\right\| _{q+1}^{q+1}\mu {(\Omega)}^{\frac
q{q+1}},
\end{eqnarray*}

where $\chi _I$ denotes the characteristic function of the
interval $I$. Inserting this estimate in (\ref{8}) we obtain

\begin{equation}\label{9}
\int_{|u| <1}|u| ^{p+1}dx\leq C\left\| u\right\| _{q+1}^{q+1}.
\end{equation}

Combining (\ref{7a}), (\ref{9}) and (\ref{5}) yields

\[
\lambda c\left( q,p\right) \left\| u\right\| _{q+1}^{q+1}\geq
c\left( \int_\Omega | \triangle u| dx\right) ^2,
\]

that is inequality (\ref{6}), then we can conclude as in Theorem
$2.4$ in \cite{BN} (see previous point).
\smallskip
\item[(iv)]
If $5>p>q>1$, then $q-p$<0 and (\ref{5}) implies directly

\[
\lambda c\left( q,p\right) \left\| u\right\| _{q+1}^{q+1}\geq c\left\|
\Delta u\right\| _1^2,
\]
where $q\in \left( 1,3\right] $; but this is inequality (\ref{6}),
so that we can conclude as before.
\smallskip
\item[(v)]
If $5>p>1>q>0$, (\ref{5}) implies

\[
\lambda c\left( q,p\right) \left\| u\right\| _2^2\geq c\left\| \Delta
u\right\| _1^2,
\]

since $q-p<0$; on the other hand, by standard elliptic estimates
(\ref{4}),

\[
\left\| \triangle u\right\| _1^2\geq c\left\| u\right\| _2^2,
\]

so that
\[
\lambda \geq \lambda _0.
\]

\item[(vi)]
Finally, if $p\geq 5>q,$ (\ref{5}) implies either

\[
\lambda c\left( q,p\right) \left\| u\right\| _{q+1}^{q+1}\geq c\left\|
\Delta u\right\| _1^2,
\]

if $q\in \left( 1,3\right] $, or

\[
\lambda c\left( q,p\right) \left\| u\right\| _2^2\geq c\left\| \Delta
u\right\| _1^2
\]

if $q\in (0,1).$ In both cases we can conclude as previously.
\end{enumerate}

The proof of Theorem \ref{t1} is now complete.

\subsection{The case $N=4$.}

We assume here that $N=4$ and

\begin{equation}\label{N=4}
g(u) =\left\{
\begin{array}{l}
|u|^{p-1}\cdot u\hspace{0.5in}|u|<1,\;\; 1<p \\
|u|^{q-1}\cdot u\hspace{0.5in}|u|\geq 1,\;\; 0<q<1
\end{array}
\right.
\end{equation}

Then we have the following result.

\begin{theorem}\label{t2}
Let $\Omega $ be strictly starshaped about the origin; suppose
that $u$ is a solution of problem (\ref{1}), with $g$ given by
(\ref{N=4}). Then
\[
\lambda \geq \lambda _0(q,p,\Omega ) >0.
\]
\end{theorem}

\emph{Proof of Theorem \ref{t2}. }By Pohozaev's identity
(\ref{pohoz}), since $\Omega $ is strictly starshaped, we have

\begin{equation}\label{10}
\begin{array}{l}
\lambda \displaystyle \frac{3-q}{q+1}\int_{|u| \geq 1}|u|
^{q+1}dx+\lambda
\frac{3-p}{p+1}\int_{|u| <1}|u|^{p+1}dx \\
\\
+4\lambda \cdot \mu \left\{ x\in \Omega :|u|\geq 1\right\}
\displaystyle \frac{q-p}{(q+1)(p+1)}\geq c\left( \int_\Omega
|\triangle u| dx\right) ^2.
\end{array}
\end{equation}

Observe that $q-p<0$ since $0<q<1<p$, from (\ref{N=4}). We discuss
separately the different cases.

\begin{enumerate}
\item[(i)]
If $1<p\leq 2$, (\ref{10}) implies

\[
\lambda \frac{3-q}{q+1}\int_{|u|\geq 1}|u|^{q+1}dx+\lambda
\frac{3-p}{p+1}\int_{|u|<1}|u|^{p+1}dx\geq c\left\| \triangle
u\right\| _1^2
\]

so that

\begin{equation}\label{11}
\lambda c\left( q,p\right) \left\| u\right\| _{q+1}^{q+1}\geq
c\left\| \triangle u\right\| _1^2,
\end{equation}

since $|u|^{p+1}\leq |u|^{q+1}$ if $|u|<1,$ and

\begin{equation}\label{12}
\lambda c\left( q,p\right) \left\| u\right\| _{p+1}^{p+1}\geq
c\left\| \triangle u\right\| _1^2
\end{equation}

since $|u|^{q+1}\leq |u|^{q+1}$ if $|u|\geq 1$. From problem
(\ref{1}) with $g$ given by (\ref{N=4}) we also obtain

\[
\left\| \triangle u\right\| _1^2\geq \left\| u\right\| _3^6
\]

while, by standard elliptic estimates on $\triangle ^{-1}$ in
$\mathbb{R} ^4$, (\ref{4}),

\[
\left\| \triangle u\right\| _1^2\geq c\left\| u\right\| _r^2
\]

for all $r\in [ 1,2)$; combining this inequality with (\ref{11})
yields

\begin{equation}\label{13}
\lambda c\left( q,p\right) \geq \left\| u\right\| _{q+1}^{1-q};
\end{equation}

(\ref{13}) implies $\left\| u\right\| _{q+1}\rightarrow 0 $ as
$\lambda \rightarrow 0$. Assume now that there exists a solution
$u$ for problem (\ref{1}); by (\ref{11}), for every $\lambda >0$
and near $0$ also $\left\| \triangle u\right\| _1^2\rightarrow 0$,
so that $\left\| u\right\| _3\rightarrow 0$ as $\lambda
\rightarrow 0$. Then, let us consider the interpolation inequality

\begin{equation}\label{14}
 \left\| u\right\| _{p+1}\leq \left\| u\right\| _r^a\cdot \left\|
u\right\| _3^{1-a}
\end{equation}

with $\displaystyle \frac{1}{p+1} =\frac{a}{r}+\frac{1-a}{3}$,
$1<r<p+1<3$; combining (\ref{14}) with the previous inequalities and %
(\ref{12}), one finds

\begin{equation}\label{15}
\left\| u\right\| _{p+1}\leq c\lambda ^{\frac{2a+1}6}\left\|
u\right\| _{p+1}^{\left( p+1\right) \left( 2a+1\right) /6}.
\end{equation}

Solving

\[
\left\{
\begin{array}{c}
\displaystyle \frac{1}{p+1}=\frac{a}{r}+\frac{1-a}{3} \\
\displaystyle \left( p+1\right) \frac{2a+1}{6}>1
\end{array}
\right.
\]

we obtain

\[
r<\frac{5-p}{3-p}\;\;\;,\;\;\; a=\frac{5-p}{2(p+1)}
\]

which satisfy $a<1,1<r<2<p+1$ if $1<p\leq 2$; inserting in
(\ref{15}) we have that there are two constants $\alpha ,\beta
>0$ such that

\begin{equation}\label{16}
c\leq \lambda ^\alpha \left\| u\right\| _{p+1}^\beta
\end{equation}

if $1<p\leq 2.$ But we have proved that $\left\| u\right\|
_3\rightarrow 0$ as $\lambda \rightarrow 0$, which implies
$\left\| u\right\| _{p+1}\rightarrow 0$ as $\lambda \rightarrow
0$, for $1<p\leq 2:$ combining this relation with (\ref{16}) we
obtain a contradiction: hence $\lambda \geq \lambda _0$.
\smallskip
\item[(ii)]
If $2<p<3$, observe that

\[
\lambda \frac{3-q}{q+1}\int_{|u|\geq 1}|u|^{q+1}dx+\lambda
\frac{3-p}{p+1}\int_{|u|<1}|u|^{p+1}dx\leq \lambda c(q,p) \left\|
u\right\| _{s+1}^{s+1}
\]

with $s\in (q,p) \cap (1,2)$; hence

\[
\lambda \left\| u\right\| _{s+1}^{s+1}\geq c\left\| \triangle
u\right\| _1^2,
\]

and we can repeat the proof given in previous point with $s$
instead of $p$, obtaining

\[
c\leq \lambda ^\alpha \left\| u\right\| _{s+1}^\beta .
\]
On the other hand,

\[
\left\| u\right\| _{s+1}^2\leq c\left\| u\right\| _3^2\leq
c\left\| \triangle u\right\| _1^2;
\]

since (\ref{11}) - (\ref{13}) are still verified, we can conclude
as in the previous point.

\smallskip
If $p\geq 3$, then (\ref{10}) implies

\[
\lambda \frac{3-q}{q+1}\int_{|u|\geq 1}|u|^{q+1}dx\geq c\left(
\int_\Omega |\triangle u|dx\right) ^2,
\]

hence (\ref{11}) - (\ref{13}) are still verified and $\left\|
u\right\| _3\rightarrow 0$ as $\lambda \rightarrow 0.$ By
(\ref{10}) we also have

\[
\lambda \left\| u\right\| _{s+1}^{s+1}\geq c\left\| \triangle
u\right\| _1^2
\]

since $|u|^{q+1}\leq |u|^{s+1}$ if $|u|\geq 1$, for all $s\in
(1,2)$; then we can conclude as in the previous points.
\end{enumerate}
The proof of Theorem \ref{t2} is now complete.

\end{document}